\documentclass[12pt]{amsart}

\usepackage{amssymb}
\usepackage{epsf}

\numberwithin{equation}{section}
\newtheorem{thm}{Theorem}[section]
\newtheorem{prop}[thm]{Proposition}
\newtheorem{lem}[thm]{Lemma}

\newtheorem{definition}[thm]{Definition}
\newenvironment{defn}{\begin{definition}\rm}{\end{definition}}

\newtheorem{rem}[thm]{Remark}
\newtheorem{example}[thm]{Example}

\newenvironment{Remark}{\begin{rem}\sl}{\end{rem}}

\newcommand{\Span}[1]{\left\langle #1\right\rangle}

\def\HP{{\mathcal H}P}
\def\Q{{\mathcal Q}\mbox{\it sym}}
\def\H{{\mathcal H}}
\def\F{K}
\def\PP{{\mathbb P}}

\def\Z{{\mathbb Z}}
\def\QQ{{\mathbb Q}}

\def\h{\overline{h}}

\newcommand{\tto}{{\relbar\joinrel\to}}
\newcommand{\ttto}{{\relbar\joinrel\relbar\joinrel\to}}

\newcommand{\QED}{
\setlength{\unitlength}{1.0pt}
\begin{picture}(7.5,7.5)
\put(5,-5){\rule{2.5pt}{2.5pt}}
\put(2.5,-2.5){\rule{5pt}{2.5pt}}
\put(2.5,0){\rule{2.5pt}{2.5pt}}
\put(0,2.5){\rule{5pt}{2.5pt}}
\end{picture}\vspace{10pt}}

\begin{document}

\title[Noncommutative Pieri operators on posets]{Noncommutative Pieri
operators on posets}

 \author[N.~Bergeron, S.~Mykytiuk, F.~Sottile,
\and S.~van Willigenburg]{Nantel Bergeron \and
Stefan Mykytiuk \and Frank Sottile \and Stephanie~van~Willigenburg}

\address[Nantel Bergeron and Stefan Mykytiuk and Stephanie van Willigenburg]
{Department of Mathematics and Statistics\\ York University\\
Toronto, Ontario M3J 1P3\\
CANADA}
\address[Frank Sottile]{Department of Mathematics and Statistics\\
University of Massachusetts\\
Amherst, MA 01003\\
USA}
\email[Nantel Bergeron]{bergeron@mathstat.yorku.ca}
\urladdr[Nantel Bergeron]{http://www.math.yorku.ca/bergeron}
\email[Stefan Mykytiuk]{mykytiuk@mathstat.yorku.ca}
\email[Frank Sottile]{sottile@math.umass.edu}
\urladdr[Frank Sottile]{http://www.math.umass.edu/\~{}sottile}
\email[Stephanie van Willigenburg]{steph@mathstat.yorku.ca}
\urladdr[Stephanie van Willigenburg]
{http://www-theory.dcs.st-and.ac.uk/\~{}stephw}

\date{19 July 2000}
\thanks{Bergeron supported in part by NSERC}
\thanks{Sottile supported in part by NSF grant DMS-9701755 and NSERC grant
        OGP0170279}
\thanks{van Willigenburg is supported in part by the Leverhulme Trust}

\thanks{J.~Comb.~Th., Ser.~A, to appear.  \copyright 2000 Academic Press}

\begin{abstract}
 We consider graded representations of the algebra $NC$ of noncommu\-tative
symmetric functions on the $\Z$-linear span of a graded poset $P$.
The matrix coefficients of such a
representation give a Hopf morphism from a Hopf algebra $\HP$ generated by
the intervals of $P$ to the Hopf algebra of quasi-symmetric functions.
This provides a unified construction of quasi-symmetric generating functions
from  different branches of algebraic combinatorics, and this
construction is useful for transferring techniques and ideas
between these branches.
In particular we show that the (Hopf) algebra of  Billera and Liu
related to Eulerian posets is dual to the peak (Hopf)
algebra of Stembridge related to enriched P-partitions, and connect this to
the combinatorics of the Schubert calculus for isotropic flag manifolds.
\end{abstract}

\maketitle

\noindent{\it Dedicated to the memory of Dr.~Gian-Carlo Rota, who inspired
us to seek the algebraic structures underlying combinatorics.}

%%%%%%%%%%%%%%%%%%%%%%%%%%%%%%%%%%%%%%%%%%%%%%%%%%
\section{Introduction}
 The algebra $\Q$ of quasi-symmetric functions was introduced by
Gessel~\cite{Ges} as a source
of generating functions for $P$-partitions~\cite{Stanley72}.
Since then, quasi-symmetric
functions have played an important role as generating
functions in combinatorics~\cite{Stanley_EC1,Stanley_ECII}.
The relation of $\Q$ to the more familiar algebra of symmetric functions was
clarified by Gelfand {\it et.~al.~}~\cite{GKal} who defined the graded Hopf
algebra $NC$ of noncommutative symmetric functions and identified $\Q$ as
its Hopf dual.

Joni and Rota~\cite{Joni_Rota} made the fundamental observation that
many discrete structures give rise to natural Hopf algebras whose coproducts
encode the disassembly of these structures (see also~\cite{Schmitt}).
A seminal link between these theories was shown by
Ehrenborg~\cite{Ehrenborg}, whose flag $f$-vector quasi-symmetric function
of a graded poset gave a Hopf morphism from a Hopf algebra of graded posets
to $\Q$.
This theory was augmented in~\cite{BS99b} where it was shown that the
quasi-symmetric function associated to an edge-labelled
poset similarly gives a Hopf morphism.
That quasi-symmetric function generalised a quasi-symmetric function
encoding the structure of the cohomology of a flag manifold as a module
over the ring of symmetric functions~\cite{BS98,BS_skew}.

We extend and unify these results by means of
a simple construction.
Given a graded representation of $NC$ on the $\Z$-linear span $\Z P$ of a
graded poset $P$, the matrix coefficients of such an action are linear maps
on $NC$ and hence quasi-symmetric functions.
In Section 2 we show how this situation gives
rise to a Hopf morphism as before.
In Section~3, we extend this
construction to an arbitrary oriented multigraph $G$.
Sections 4, 6, and 7 give examples of this construction, including
rank selection in posets, flag $f$-vectors of polytopes,
$P$-partitions, Stanley symmetric functions, and the multiplication of
Schubert classes in the cohomology of flag manifolds.

In Section~5, we discuss how properties
of the combinatorial structure of $G$ may be understood through the
resulting quasi-symmetric function.
This analysis allows us to relate
work of Bayer, Billera, and Liu~\cite{BaBi,BilLiu} on Eulerian posets
with work of Stembridge~\cite{Stembridge_enriched} on enriched
P-partitions. More precisely, we show that the quotient of $NC$
by the ideal of the generalised Dehn-Somerville relations is dual to the
Hopf subalgebra of peak functions in $\Q$.
We also solve the conjecture of \cite{BMSW}, showing that the shifted
quasi-symmetric
functions form a Hopf algebra.
These functions were introduced by Billey and
Haiman~\cite{BilHaim} to define Schubert polynomials for all types.

In Section 7, we show how a natural generating function for enumerating peaks
in a labelled poset is the quasi-symmetric function for an enriched
structure on the poset.
Special cases of this combinatorics of peaks include Stembridge's theory of
enriched $P$-partitions~\cite{Stembridge_enriched}, the Pieri-type formula
for type $B$ and $C$ Schubert polynomials in~\cite{BS_lag-pieri}, and
Stanley symmetric functions of types $B$, $C$, and $D$.
These examples linking the diverse areas of Schubert calculus,
combinatorics of polytopes and $P$-partitions
illustrate how this theory  transfers techniques and ideas
between disparate areas of combinatorics.

We thank Sarah Witherspoon who contributed to the appendix on Hopf
algebras, and Geanina Tudose for her
assistance with fusion coefficients.

%%%%%%%%%%%%%%%%%%%%%%%%%
\section{Pieri operators on posets}

Many interesting families of combinatorial constants can be understood as an
enumeration of paths in a ranked partially ordered set (poset) which satisfy
certain conditions.
One example of this is the  Littlewood-Richardson rule in the theory of
symmetric functions \cite{Mac}.
This rule describes the multiplicity $c_{\mu,\nu}^{\lambda}$ of a Schur
function $S_\lambda$ in the product $S_\mu S_\nu$ of two others.
The constants $c_{\mu,\nu}^{\lambda}$ can be seen as an enumeration of all
paths in Young's lattice from $\mu$ to $\lambda$ satisfying some conditions
imposed by $\nu$.
We note that the constants $c_{\mu,\nu}^{\lambda}$ are invariant under
certain isomorphisms of intervals in Young's lattice,
namely $c_{\mu,\nu}^{\lambda} = c_{\tau,\nu}^{\pi}$ whenever
$\lambda/\mu = \pi/\tau$.
The skew Schur functions $S_{\lambda/\mu}$ are generating
functions of these constants as we have
 $$
  S_{\lambda/\mu}\ =\ \sum_\nu c_{\mu,\nu}^{\lambda}\, S_\nu .
 $$

We generalise the principles of this example,
introducing families of algebraic operators to {\sl select} paths in a given
poset.
Here, analogues of the Littlewood-Richardson constants count paths in
the poset satisfying some conditions imposed by the family of operators.
These enumerative combinatorial invariants of
the poset are encoded by
generating functions which generalise the skew Schur functions.
We show that the association of such a generating function to a
poset induces a Hopf morphism to $\Q$.

\smallskip
Let $(P,<)$ be a graded
poset with rank function $rk\colon P\to\Z^+$ and let $\Z P$ be the free
graded $\Z$-module generated by the elements of $P$.
For an integer $k>0$, a  (right) {\sl Pieri operator} on $P$
is a linear map $\h_k\colon\Z P\to\Z P$ which respects the poset structure.
By this we mean that for all $x\in P$, the support
of $x.\h_k\in \Z P$ consists only of elements $y\in P$ such that $x<y$ and
$rk(y)-rk(x)=k$.
We note that such an operator $\h_k$ is of degree $k$ on $\Z P$.

Gelfand {\it et.~al.~}\cite{GKal} define the Hopf algebra $NC$ of
noncommutative
symmetric functions to be the free associative algebra
$\Z\langle h_1,h_2,\ldots\rangle$ with
a generator $h_k$ in each positive degree $k$ and coproduct
$\Delta h_k=\sum_{i=0}^k h_i\otimes
h_{k-i}$, where $h_0=1$.
It follows that given a family of Pieri operators $\{\h_k\}_{k> 0}$
on a poset $P$, the map $h_k\mapsto \h_k$ turns $\Z P$ into
a graded (right) $NC$-module.
Conversely, any graded right action of $NC$ on $\Z P$ which
respects the poset structure of $P$ gives a family of Pieri operators on
$P$.
When the context is clear, we may identify the generator $h_k$ with the
operator $\h_k$.

 Given such a representation of $NC$ on $\Z P$ and $x,y\in P$, the
association of $\Psi\in NC$
to the coefficient of $y$ in $x.\Psi$ is a linear map on $NC$.
These matrix coefficients are elements of the Hopf dual of $NC$ which is the
Hopf algebra $\Q$ of quasi-symmetric functions~\cite{GKal}.
These coefficients vanish unless $x\leq y$.
This gives a collection
of quasi-symmetric functions $\F_{[x,y]}$ associated to every interval
$[x,y]$ of $P$.

 Let $\HP$ be the free $\Z$-module with basis given by Cartesian products of
intervals $[x,y]$ of $P$, modulo identifying all singleton intervals $[x,x]$
with the unit $1$ and empty intervals with zero.
Then $\HP$ is a graded $\Z$-algebra whose product is the cartesian
product of intervals, and whose grading is induced by the rank of an
interval of $P$.
It has a natural coalgebra structure induced by
 $$
  \Delta A\ =\ \sum_{x\in A}\, [\hat{0}_A,x]\otimes[x,\hat{1}_A],
 $$
where $A=[\hat{0}_A,\hat{1}_A]$ is an interval of $P$ with minimal element
$\hat{0}_A$ and maximal element $\hat{1}_A$.
Projection onto $\Z$ of the degree 0 component of $\HP$ is the
counit. It follows that $\HP$ is a bialgebra. It is graded, therefore by
Proposition~\ref{GradedHopf}, there is a unique antipode and $\HP$ is a Hopf
algebra.

  \begin{thm}
    For any graded poset $P$, $\HP$ is a Hopf algebra.
  \end{thm}

Suppose we have a family of Pieri operators on a poset $P$.
Since $NC$ is a Hopf algebra, the action of $NC\otimes NC$ on
$\Z P\otimes\Z P= \Z(P\times P)$ pulls back along the the coproduct $\Delta$
to give an action of $NC$ on $\Z(P\times P)$.
We iterate this and use coassociativity to get an action of $NC$ on the
$\Z$-linear span of  $P^k$.
Since a product of intervals of $P$ is an interval in such an iterated
product of $P$
with itself, we may extend the definition of $K$ to the generators of $\HP$
and then by linearity to $\HP$ itself, obtaining a $\Z$-linear homogeneous
map $\F\colon\HP\to \Q$.
Let $\Span{\cdot,\cdot}$ be the bilinear form on $\Z P$ induced by the
Kronecker delta function on the elements of $P$.

\begin{thm}\label{Hopfmap}
The map $\F\colon\HP\to \Q$ is a morphism of Hopf algebras.
\end{thm}

\proof
 We show that $\F$ respects product and coproduct, which suffices.
For $x,y\in P$ and $\Psi\in NC$, we have
$K_{[x,y]}(\Psi)= \Span{x.\Psi,\, y}$.
Thus for $x\in P$ and $\Psi\in NC$,
 $$
  x.\Psi\ =\ \sum_y \Span{x.\Psi,y}\,  y\ =\
      \sum_y \F_{[x,y]}(\Psi)\, y.
 $$

Let $A=[\hat{0}_A,\hat{1}_A]$ and $B=[\hat{0}_B,\hat{1}_B]$ be intervals of
$P^{k_1}$ and $P^{k_2}$ respectively. For $\Psi\in NC$,
using Sweedler notation for the coproduct
 $$
  \Delta \Psi\ =\ \sum \Psi_a\otimes \Psi_b,
 $$
and the duality between the product of $\Q$ and the coproduct of $NC$, we
obtain
 \begin{eqnarray*}
   \F_{A\times B}(\Psi)
   &=& \Span{(\hat{0}_A\otimes\hat{0}_B).\Psi,\ \hat{1}_A\otimes\hat{1}_B}\\
   &=& \Span{\sum\hat{0}_A.\Psi_a\otimes\hat{0}_B.\Psi_b,\
   \hat{1}_A\otimes\hat{1}_B}\\
   &=& \sum \Span{\hat{0}_A.\Psi_a,\ \hat{1}_A}\Span{\hat{0}_B.\Psi_b,\
   \hat{1}_B}\\
   &=& \sum \F_A(\Psi_a)\F_B(\Psi_b)
   \ =\ (\F_A\otimes \F_B)(\Delta\Psi)
   \ =\ (\F_A\cdot \F_B)(\Psi).
 \end{eqnarray*}

 Let $A=[\hat{0}_A,\hat{1}_A]$ be an interval of $P^k$ and $\Psi,\Phi\in
NC$.
Using the duality
between the coproduct of $\Q$ and the product of $NC$, we have
 \begin{eqnarray*}
   (\Delta \F_A)(\Psi\otimes \Phi)
   &=& \F_A(\Psi\cdot\Phi)
      \ =\ \Span{(\hat{0}_A.\Psi).\Phi,\ \hat{1}_A}\\
   &=& \Span{\sum_y(\F_{[\hat{0}_A,y]}(\Psi) y).\Phi,\ \hat{1}_A}\\
   &=& \Span{\sum_{x,y}\F_{[\hat{0}_A,y]}(\Psi)\F_{[y,x]}(\Phi) x,\
        \hat{1}_A}\\
   &=& \sum_{y\in A}\F_{[\hat{0}_A,y]}(\Psi)\F_{[y,\hat{1}_A]}(\Phi)
      \ =\ \F_{\Delta A}(\Psi\otimes \Phi).
      \qquad \QED
 \end{eqnarray*}

The map $K$ is a generating function for the enumerative combinatorial
invariants associated to the $NC$-structure of $\Z P$.
Let $\{ a_\alpha\}$ be a graded basis of $NC$ and let $\{ b_\alpha\}$
be the corresponding dual basis in $\Q$.  Then
\begin{equation}\label{eq:Kformula}
\F_{[x,y]}\ =\ \sum_\alpha\  \langle x.\overline{a}_\alpha,y \rangle\,
b_\alpha.
\end{equation}
We interpret the coefficient of $b_\alpha$ in $\F_{[x,y]}$
as the number of paths from $x$ to $y$ satisfying some condition imposed by
$a_\alpha$.

We reformulate Equation~\ref{eq:Kformula} in terms of the Cauchy
element of Gelfand {\it et.~al.~}~\cite{GKal}.
This element relates each graded basis of the Hopf algebra
$NC=\bigoplus_{n\ge 0} NC_n$ to its corresponding
dual basis in the Hopf algebra $\Q=\bigoplus_{n\ge 0} \Q_n$.
More precisely, let $\alpha=(\alpha_1,\alpha_2,\ldots,\alpha_\ell)$ with
$\ell\ge 0$  be a sequence of positive integers.
Such a sequence is a {\sl composition} of $n$, denoted
$\alpha\models n$, if $n=\sum_{i=1}^\ell \alpha_i$. By convention the empty
sequence for
$\ell=0$ is the unique composition of $0$.
The complete $NC$-functions
$\{S^\alpha=h_{\alpha_1}h_{\alpha_2}\cdots h_{\alpha_\ell}\}_{\alpha\models
n}$ and the ribbon $NC$-functions $\{R_\alpha\}_{\alpha\models n}$ form two
bases of $NC_n$.
Similarly the  monomial quasi-symmetric functions
$\{M_\alpha\}_{\alpha\models n}$ and the complete
quasi-symmetric functions $\{F_\alpha\}_{\alpha\models n}$ form two bases of
$\Q_n$.
In the graded completion of
$\bigoplus_{n\ge 0} NC_n\otimes {\Q}_n$ we have the following {\sl Cauchy
element}:
  \begin{equation}\label{eq:NCcauchy}
  {\mathcal C}\ :=\ \sum_\alpha a_\alpha\otimes b_\alpha\
                 =\ \sum_\alpha R_\alpha\otimes F_{\alpha}\
                 =\ \sum_\alpha S^{\alpha}\otimes M_\alpha,
  \end{equation}
where $\{a_\alpha\}$ and $\{b_\alpha\}$ is any pair of dual graded bases.

The right action of $NC$ on $\Z P$ extends linearly to an
action of the graded completion of
$\bigoplus_{n\ge 0} NC_n\otimes {\Q}_n$ on the completion of $\Q
\otimes_{\Z} \Z P$.
The following theorem is simply a reformulation of Equation~\ref{eq:Kformula}.

\begin{thm}\label{Cauchy}
For any family of Pieri operators on a poset $P$ and $x,y\in P$, $$
\F_{[x,y]}\ =\ \Span{x.{\mathcal C},\ y}. $$
\end{thm}

Expanding the quasi-symmetric function $\F_{[x,y]}$ in a basis of $\Q$ gives
a family of enumerative combinatorial invariants for the given action
of $NC$ on the poset $P$.
In this way, the functions $\F_{[x,y]}$ are seen to be analogues of the
skew Schur functions presented at the beginning of this section.

%%%%%%%%%%%%%%%%%%%%%%%%%
\section{Pieri operators on graphs}

We extend this simple construction on graded posets to
(locally finite) oriented multigraphs.
Let $G=(V,E)$ be a multigraph where $V$ is the set of vertices
and $E$ is a function $V\times V\to\Z^+$ such that
$$
  \sum_{y'\in V}E(x,y') \quad\hbox{and}\quad \sum_{x'\in V}E(x',y)
$$
are both finite for all
$x,y\in V$.
The value $E(x,y)$ identifies the number of arrows from $x$ to $y$ in $G$.
The function $E$ is the incidence matrix of the graph $G$, and
$E^r$ is the matrix product of $r$ copies of $E$.
Given $x,y\in V$, let $[x,y]$ be the set
of all paths from $x$ to $y$.
Consider a graded version of this set,
$$
  [x,y]\ =\ \bigcup_{r\ge 0} [x,y]^{(r)},
$$
where the interval $[x,y]^{(r)}$ is the set of all paths of
length $r$ from $x$ to $y$.
Note that $|[x,y]^{(r)}|=E^r(x,y)$ is finite.

Let $\Z G$ denote the free $\Z$-module generated by $V$.
Here, a {\sl Pieri operator} is a linear map
$\h_k\colon \Z G\to \Z G$ where for all $x\in V$, the support of
$x.\h_k\in\Z G$
consists of elements $y\in V$ such that $E^k(x,y)>0$.
As before, a family $\{\h_k\}_{k> 0}$ of
Pieri operators induces on  $\Z G$ the structure of an $NC$-module.
We thus obtain a collection of
linear maps $\F_{[x,y]^{(r)}}\colon NC\to \Z$ given by $\Psi\mapsto \langle
x.\Psi^{(r)},y\rangle$ where $\Psi^{(r)}$ is the $r$th-homogeneous component
of $\Psi$, and thus  quasi-symmetric functions $\F_{[x,y]^{(r)}}\in\Q_r$.

We define a Hopf algebra $\H G$ associated to $G$.
Define the product of two intervals by
$[x,y]^{(r)}\times [u,v]^{(s)}:=[(x,u),(y,v)]^{(r+s)}$,
an interval in $G\times G$.
Let $\H G$ be the free $\Z$-module with basis given by products of intervals
$[x,y]^{(r)}$ in $G$, modulo identifying all intervals
$[x,x]^{(0)}$ with the unit $1$
and all empty intervals $[x,y]^{(r)}$ with zero.
If we let $r$ be the degree of an element
$[x,y]^{(r)}$, then $\H G$ is a graded $\Z$-algebra with product $\times$.
The algebra $\H G$ has a natural coalgebra structure induced by
$$
  \Delta [x,y]^{(r)}
  = \sum_{s=0}^r \sum_{z\in V} [x,z]^{(s)}\otimes [z,y]^{(r-s)}.
$$
The counit is again the projection onto $\Z$ of the degree $0$ component.
Since the bialgebra $\H G$ is graded, we have the following theorem.

   \begin{thm}
    For any oriented multigraph $G$, $\H G$ is a Hopf algebra.
   \end{thm}

Suppose we have a family of Pieri operators on a graph $G$.
As in Section 2, we have an action of $NC$ on
the $\Z$-linear span of $G^k$, for any positive integer $k$.
Since the generators of $\H G$ are sets of the form $[w,z]^{(t)}$
in $G^k$, we may define the quasi-symmetric function $K$ on each of these
generators of $\H G$, and then extend by linearity to $\H G$ itself,
obtaining a $\Z$-linear graded map $K:\H G \to \Q$.
We leave to the reader the straightforward extension of
Theorem~\ref{Hopfmap}.

\begin{thm}\label{GHopfmap}
The map $\F\colon\H G\to \Q$ is a morphism of Hopf algebras.
\end{thm}

To extend Theorem~\ref{Cauchy}, we decompose the Cauchy element
${\mathcal C}$ into
its homogeneous components: ${\mathcal C}=\sum_{r\ge 0} {\mathcal C}_r$.
For example, we can
use ${\mathcal C}_r=\sum_{\alpha\models r} S^{\alpha}\otimes M_\alpha$.
The following is
immediate.

  \begin{thm}\label{GCauchy}
    For any family of Pieri operators on a graph $G$ and $x,y\in G$,
    $$
     \F_{[x,y]^{(r)}}\ =\ \Span{x.{\mathcal C}_r,\ y}. $$
  \end{thm}

\begin{rem}
\rm These constructions generalise those of
Section~2.
Given a ranked poset $P$ we  associate to it the incidence graph
$G_P=(P,E)$ where $E(x,y)=1$ if $y$ covers $x$ and
$E(x,y)=0$ otherwise.
The intervals $[x,y]^{(r)}$ are empty unless
$r=rk(y)-rk(x)$ in which case
$[x,y]^{(r)}$ is equal to the saturated chains in $[x,y]$.
In such a case we omit the superscript $(r)$ and arrive again at the results
of Section 2.
\end{rem}

%%%%%%%%%%%%%%%%%%%%%%%%%
\section{Three simple examples}

We give three simple examples to illustrate our theory.

\begin{example}\label{Ex:SPE}
Simple path enumeration. \quad\rm
Given a graph $G$, define the Pieri operator $\h_k\colon\Z G\to\Z G$ by
$$
x.\h_k = \sum _{y\in G} E^k(x,y)y.
$$
This action of $NC$ satisfies $\h_a\h_b=\h_{a+b}$ for all $a,b\in \Z^+$.
{From} this we deduce that
$\F_{[x,y]^{(r)}}= E^r(x,y)\sum _{\alpha\vDash r} M_\alpha$.
Thus $\F$ simply enumerates  all paths of length $r$ from $x$ to $y$.
When $G$ is a graded poset $P$, $E^r(x,y)=0$ unless $r=rk(y)-rk(x)$ and
$x\leq y$ in $P$.
In this case, $E^r(x,y)$ counts the saturated chains in the interval $[x,y]$.
\end{example}

\begin{example}\label{Ex:SkewSchur} Skew Schur functions. \quad\rm
Let $(P,<)$ be Young's lattice of partitions. For $\mu\in P$,
define $\mu.\h_k$ to be the sum of all partitions $\lambda$ such that
$\lambda/\mu$
is a horizontal strip and $|\lambda|-|\mu|=k$.
This  family of Pieri operators lifts the action of
the algebra $\Lambda$ of symmetric functions on itself.
It follows that $\F_{[\mu,\lambda]}$ is the skew Schur function
$S_{\lambda/\mu}$.
\end{example}

\begin{example}\label{Ex:Flagf} Rank selection Pieri operators and flag
$f$-vectors. \quad\rm
Given any ranked poset $P$, consider the Pieri operator obtained by setting
$x.\h_k$ equal to the
sum of all $y>x$ such that $rk(y)-rk(x)=k$.
In this case, $\langle x.\overline{S}^\alpha, y\rangle$ counts all chains
in the rank-selected poset obtained from $[x,y]$ with ranks given by
$\alpha$, and $\F$ is Ehrenborg's flag
$f$-vector quasi-symmetric generating function~\cite{Ehrenborg}.
\end{example}

%%%%%%%%%%%%%%%%%%%%%%%%%
\section{Structure from Hopf subalgebras}
Suppose, for a family of posets, we have a class of enumerative
combinatorial invariants
 which possesses some additional structure.
In many situations, the associated families of Pieri operators satisfy some
relations, and the resulting actions of $NC$ are carried by a Hopf
quotient of $NC$.
Equivalently, the images of $K$ lie in the dual of this quotient, a Hopf
subalgebra of $\Q$.

More precisely, let
an action of $NC$ on $\Z G$ be given by a homomorphism $\phi$ from $NC$ to
the linear endomorphism ring End$(\Z G)$.
Let $\mathcal I$ be an ideal generated by some relations
satisfied by the Pieri operators.
When ${\mathcal I}$ is a Hopf ideal, so that we have
$\Delta({\mathcal I}) \subset
   {\mathcal I}\otimes NC + NC\otimes{\mathcal I}$,
we have the commuting diagram
$$
\begin{picture}(107,58)
\put(10,45){$NC$}	\put(85,45){End$(\Z G)$}
\put(6,0){$NC$\big/\lower4pt\hbox{$\mathcal I$}}
\put(45,51){\scriptsize $\phi$}
\put(30,48){\vector(1,0){50}}  \put(23,40){\vector(0,-1){27}}
 \put(40,10){\vector(3,2){46}} \put(51,28){\scriptsize $\phi^*$}
\end{picture}
\hskip50pt
\raise28pt\hbox{$\implies$}
\hskip20pt
\begin{picture}(107,58)
\put(10,45){$\H G$}	\put(85,45){$\Q$}
\put(85,0){$\left(\hbox{$NC$\big/\lower4pt\hbox{$\mathcal I$}}\right)^*$}

\put(45,51){\scriptsize $\F$}
\put(32,48){\vector(1,0){50}}  \put(100,23){$\cup$}
 \put(32,40){\vector(3,-2){46}} \put(41,18){\scriptsize $\F$}
\end{picture}
$$
as the functions $K_{[w,z]^{(t)}}$ are characters of representations
$\phi^{\otimes k}$ on $\Z G^k$.

In particular, Equation~\ref{eq:Kformula} has a more specialised form.
Given a basis
$\{c_\lambda\}$ of $NC\big/\lower4pt\hbox{$\mathcal I$}$ and $\{d_\lambda\}$
its dual basis inside $\Q$, we have
  \begin{equation}\label{eq:subKformula}
   \F_{[x,y]^{(r)}}=\sum_\lambda \langle x.c_\lambda,y\rangle  d_\lambda
  \end{equation}
where the sum is only over the index set of the given basis for
$NC\big/\lower4pt\hbox{$\mathcal I$}$.
Here the numbers $\langle x.c_\lambda,y\rangle$ are special cases of the
enumerative invariants in Equation~\ref{eq:Kformula}.

We illustrate these principles in a series of examples which introduce
certain classes of Pieri operators defined by quotients of
$NC$.

\begin{example} Simple path enumeration.\rm\quad
In Example~\ref{Ex:SPE}, the ideal
${\mathcal I}$ is generated by $h_{a+b}-h_ah_b$ for all $a,b>0$.
This is not a Hopf ideal since
  \begin{eqnarray*}
    \Delta(h_2-h_1h_1)
         &=&    h_2\otimes 1 + h_1\otimes h_1 + 1\otimes h_2
              -(h_1\otimes 1 + 1\otimes h_1)^2\\
         &=&  (h_2-h_1h_1)\otimes 1 + 1\otimes (h_2-h_1h_1)
              -h_1\otimes h_1,
  \end{eqnarray*}
which is not contained in
${\mathcal I}\otimes NC + NC\otimes{\mathcal I}$.
\end{example}

\begin{example}\label{ex:sym} Symmetric Pieri operators.\rm\quad
A family of Pieri operators is {\sl symmetric} if
$\h_a\h_b= \h_b\h_a$ for all $a,b>0$.
In this case,
$NC\big/\lower4pt\hbox{$\mathcal I$}\cong\Z[h_1,h_2,\ldots]$,
which is the self-dual Hopf algebra $\Lambda$ of
symmetric functions (see~\cite{Mac, Stanley_ECII}),
and thus ${\mathcal I}$ is a Hopf ideal.
Symmetric Pieri operators satisfy
$x.S^\alpha =x.S^\beta$ whenever $\alpha$ and $\beta$ determine the same
partition, and hence by Equation~\ref{eq:Kformula}
we can write
$\F _{[x,y]}$ in the form
  $$
   \sum _{\lambda\vdash r} A^\lambda\,
     \sum_{\lambda(\alpha)=\lambda} M_\alpha
  $$
where $r$ is the rank of the interval $[x,y]$,
$\lambda(\alpha)$ is the partition determined by $\alpha$,
and $A^\lambda$ is some constant.
By definition, $\sum _{\lambda(\alpha)=\lambda} M_\alpha$ is the symmetric
function $m_\lambda$, and so we see again that the image of $\F$
lies in $\Lambda$.
Symmetric Pieri operators can be found in Example~\ref{Ex:SkewSchur}
and in Sections 6 and 7.

It is interesting to use other known dual bases of
$\Lambda$ in Equation~\ref{eq:subKformula}, in particular, its
self-dual basis $\{S_\lambda\}$ of Schur functions.
\end{example}

\begin{example} Flag $f$-vectors of Eulerian posets.\rm\quad
Consider Example~\ref{Ex:Flagf} when the given ranked poset
$P$ is Eulerian.
The flag $f$-vectors of Eulerian posets satisfy the {\sl linear}
generalised Dehn-Sommerville or
Bayer-Billera relations~\cite{BaBi}.
Billera and Liu \cite[Proposition 3.3]{BilLiu} show that the ideal of
relations satisfied by such Pieri operators is generated by the
(even) Euler relations
  \begin{equation}\label{eq:Euler}
      \sum_{i+j=2n} (-1)^i\, \h_i\h_j \ =\  2 \h_{2n} + \sum_{i=1}^{2n-1}
      (-1)^i\,\h_i\h_{2n-i}\ = 0\,,
  \end{equation}
where $n$ is a positive integer.
As in \cite{BilLiu}, let
${\mathcal I}$ be the ideal of  $NC$ generated by
 $$
  X_{2n} \ :=\ \sum_{i+j=2n} (-1)^i\, h_ih_j\ =\
  2 h_{2n} + \sum_{i=1}^{2n-1}  (-1)^i\, h_ih_{2n-i}.
 $$
Then we have the following algebra isomorphism, 
 $$
  \QQ \otimes NC\big/\lower4pt\hbox{${\mathcal I}$}
   \ \cong\  \QQ\langle y_1,y_3,y_5,\ldots\rangle.
 $$
where $y_i$ has degree $i$.
We identify the dual
$\left(NC\big/\lower4pt\hbox{${\mathcal I}$}\right)^*$ inside
$\QQ\otimes\Q$
to be the peak Hopf algebra
$\Pi$ introduced by Stembridge \cite{Stembridge_enriched} in
his study of enriched P-partitions.
This shows that ${\mathcal I}$ is a Hopf ideal over $\QQ$.
\end{example}

   \begin{thm}\label{thm:PI}
    ${\displaystyle
     \left(\QQ\otimes NC\big/\lower4pt\hbox{${\mathcal I}$}\right)^* \
     \  \cong \ \QQ\otimes\Pi}$.
   \end{thm}

Let us clarify some notation for the proof of Theorem~\ref{thm:PI}.
Given compositions $\alpha,\beta\models m$, write
$\beta\preccurlyeq\alpha$ if $\beta$
is a refinement of $\alpha$ and let $\beta^*$
be the refinement of $\beta$ obtained by replacing all components
$\beta _i >1$ of $\beta$ for $i>1$ with $[1, \beta_i -1]$.
Given a composition $\alpha\models m$ with $\alpha_1>1$ if $m>1$,
the Billey-Haiman shifted quasi-symmetric functions~\cite{BilHaim} are
shown~\cite{BMSW}  to have the formula
  \begin{equation}\label{eq:peak}
    \theta _\alpha \ =\
     \sum_{\beta \models m \atop \beta^\ast \preccurlyeq\alpha}
      2^{k(\beta)}  M_\beta,
  \end{equation}
where  $k(\beta)$ is the number of  components of $\beta$.

If $\alpha$ is a composition with all components greater than 1, except
perhaps the last, then we call $\alpha$ a {\sl peak composition} and
$\theta_\alpha$ a {\sl peak function}.
In \cite{Stembridge_enriched} Stembridge  shows that the linear span $\Pi$
of the peak functions is a subalgebra of $\Q$.
In fact, $\Pi$ is a Hopf subalgebra of $\Q$~\cite{BMSW}.

Recall that in the identification of $\Q$ as the graded linear dual of $NC$,
the families $\{M_\alpha\}$ and $\{S^\alpha\}$ are dual bases.
That is, $M_\alpha(S^\beta)=1$ if $\alpha=\beta$ and $0$ otherwise.
Given any two compositions $\eta=(\eta_1,\eta_2,\ldots)$ and
$\epsilon=(\epsilon_1,\epsilon_2,\ldots)$, let $\eta\cdot\epsilon$
be the concatenation  $(\eta_1,\eta_2,\ldots,\epsilon_1,\epsilon_2,\ldots)$.

  \begin{lem}\label{zero-map}
    The peak algebra $\Pi$ annihilates the ideal ${\mathcal I}$.
  \end{lem}

\proof
We show that a peak function $\theta _\alpha$ annihilates any function of
the form $S^\beta X_{2n}S^\gamma$.
Since $M_\eta (S^\epsilon)=0$  unless
$\eta =\epsilon$, it follows that we need only study those  summands
$2^{k(\delta)}M_\delta$ in $\theta _\alpha$ such that
either $\delta=\beta\cdot 2n\cdot\gamma$  or else
$\delta=\beta\cdot i\cdot (2n-i) \cdot\gamma$.
Now if $2^{k(\beta\cdot 2n\cdot\gamma)}M_{\beta\cdot 2n\cdot\gamma}$
is a summand of $\theta  _\alpha$, then it follows that all summands of the
form
$2^{k(\beta\cdot i\cdot (2n-i) \cdot\gamma)}
     M_{\beta\cdot i\cdot (2n-i)\cdot\gamma}$
will also belong  to $\theta _\alpha$.
By the Euler relations~\ref{eq:Euler}, it follows immediately
that $ S^\beta X_{2n}  S^\gamma$ is annihilated by $\theta_\alpha$.

We are left to consider the case where $2^{k(\beta\cdot i\cdot (2n-i)
\cdot\gamma)}M_{\beta\cdot i\cdot (2n-i)\cdot\gamma}$ is
a summand of $\theta _\alpha$ but not $2^{k(\beta\cdot 2n\cdot\gamma)}
M_{\beta\cdot  2n\cdot\gamma}$.
Observe that from the definition~\ref{eq:peak} we must have $n>1$ since if
$2^{k(\beta\cdot 1\cdot 1\cdot\gamma)}M_{\beta\cdot 1\cdot
1\cdot\gamma}$ is a summand of $\theta _\alpha$ then
$2^{k(\beta\cdot 2\cdot\gamma)}M_{\beta\cdot 2\cdot\gamma}$
will be too as $(\beta\cdot 1\cdot
1\cdot\gamma)^* = (\beta\cdot 2\cdot\gamma)^*$.
Suppose $\beta\models m$ and let $j$ be such that
 $$
   \alpha_1+\alpha_2+\cdots+\alpha_{j-1}\le
   m<\alpha_1+\alpha_2+\cdots+\alpha_j.
 $$
Then $M_{\beta\cdot i\cdot (2n-i)\cdot\gamma}$ is in the support of
$\theta _\alpha$
if and only if $\alpha_1+\alpha_2+\cdots+\alpha_j$ is $m+i$ or $m+i+1$.
If it is $m+i$, then for $i\neq 1$
 $$
  \theta _\alpha (S^\beta X_{2n} S^\gamma)\ =\
  \theta _\alpha((-1)^i S^{\beta\cdot i\cdot(2n-i)\cdot\gamma}
       +(-1)^{(i-1)}S^{\beta\cdot i-1\cdot (2n-i+1)\cdot\gamma})
  \ =\  0.
 $$
If $i=1$, then
$(\beta\cdot 1\cdot(2n-1)\cdot\gamma)^*\preccurlyeq\alpha$.
Since $\theta_\alpha$ is a peak function and $n>1$,
we must have $\alpha_{j+1}>1$.
This implies that
$2n-1\le\alpha_{j+1}$, hence
$(\beta\cdot 2n\cdot\gamma)^*\preccurlyeq\alpha$ and
$2^{k(\beta\cdot 2n\cdot\gamma)}M_{\beta\cdot 2n\cdot\gamma}$
is a summand of $\theta_\alpha$, which contradicts our assumption.

A similar argument for $m+i+1$ completes the proof of the lemma.
\QED

\noindent{\it Proof of Theorem~\ref{thm:PI}.}
By Lemma~\ref{zero-map},
$\QQ\otimes\Pi\subseteq (\QQ\otimes NC/{\mathcal I})^* \ \cong\
  \big(\QQ\langle y_1,y_3,y_5,\ldots\rangle\big)^*$.
This containment is an equality since the dimension of the $i$th
homogeneous component of
both  $\Pi$~\cite{Stembridge_enriched} and
$ \QQ\langle y_1,y_3,y_5,\ldots\rangle$~\cite{BilLiu} is the
$i$th Fibonacci  number.
\QED

\begin{defn} \rm
Pieri operators are symmetric if the image of
$\F$ lies within the algebra $\Lambda$ of symmetric functions.
Similarly, Pieri operators are {\sl Eulerian} if the image of
$\F$ lies within $\Pi_{\mathbb Q}$, the ${\mathbb Q}$-span of $\Pi$.
This occurs if there is some scalar
multiple $\alpha_k\h_k$ of each Pieri operator such that the
$\alpha_k\h_k$ satisfy the Euler relations~\ref{eq:Euler}.
\end{defn}

We solve the conjecture presented in \cite{BMSW}
related to the general functions $\theta_\alpha$ introduced by Billey and
Haiman \cite{BilHaim}.
Let $\Xi$ be the $\QQ$-linear span of all the $\theta_\alpha$.

  \begin{thm}
   The space $\Xi$ is a Hopf subalgebra of $\QQ\otimes\Q$.
   Moreover the set
   $$
     {\mathcal J}\ =\
      \{\Psi\in NC\, |\,\, \theta(\Psi)=0 \hbox{ for all }\theta\in\Xi\}
   $$
   is the principal ideal generated by $X_2=2h_2-h_1h_1$.
  \end{thm}

\proof
We first show that ${\mathcal J}$ is an ideal and it is included in
${\mathcal I}=\langle X_{2n}\rangle$, the ideal generated by the Euler
relations.
By Theorem 3.2 of~\cite{BMSW},  $\Xi$ is a coalgebra.
Hence $\Xi^*=NC\big/\lower4pt\hbox{$\mathcal J$}$ is an algebra,
which shows that   ${\mathcal J}$ is an ideal.
Since $\Pi\subset\Xi$ we have that
${\mathcal J}\subset {\mathcal I}$.
Now it is straightforward to check that $X_2\in\mathcal J$.
Let $\hat{\mathcal J}\subseteq{\mathcal J}$ be the principal ideal
generated by $X_2$.
Since $\Delta(X_2)=1\otimes X_2 + X_2\otimes 1$ we have that
$\hat {\mathcal J}$ is a Hopf ideal and
$NC\big/\lower4pt\hbox{$\hat{\mathcal J}$}$ is a Hopf algebra.
Its dual $\left(NC\big/\lower4pt\hbox{$\hat{\mathcal J}$}\right)^*$
is a Hopf subalgebra of $\Q$ contained in $\Xi$.
To conclude our argument, we   show that the dimension of the
homogeneous components of degree $n$ in
$NC\big/\lower4pt\hbox{$\hat{\mathcal J}$}$ and $\Xi$ are equal
for all $n$.
In $NC\big/\lower4pt\hbox{$\hat{\mathcal J}$}$,
 the homogeneous component of degree $n$
has dimension given by the number of compositions of $n$ that contain no
component equal to $2$.
This satisfies the recurrence $\pi_n=\pi_{n-1}+\pi_{n-2}+\pi_{n-4}$ with
initial conditions
$\pi_1=1$, $\pi_2=1$, $\pi_3=2$ and $\pi_4=4$.
This is exactly the recurrence of Theorem 4.3 in~\cite{BMSW}
given for calculating the dimension of the homogeneous component of degree
$n$ in $\Xi$.
Hence
$\left(\QQ\otimes NC\big/\lower4pt\hbox{$\hat{\mathcal J}$}\right)^*=\Xi$
is a Hopf algebra and
$\hat{\mathcal J}={\mathcal J}$.
\QED

%%%%%%%%%%%%%%%%%%%%%%%%%
\section{Descent Pieri operators}

\begin{defn}\label{Ex:EL}
An (edge)-labelled poset is a graded poset $P$ whose
covers (edges of its Hasse diagram) are labelled with integers.
To enumerate chains according to the descents in their sequence
of (edge) labels,
we use the {\sl descent Pieri operator}
 $$
   x.\h_k\ :=\ \sum_\omega {\rm end}(\omega),
 $$
where the sum is over all chains $\omega$ of length $k$ starting at $x$,
 $$
  \omega \ :\  x  \ \stackrel{b_1}{\tto}\
               x_1\ \stackrel{b_2}{\tto}\  \cdots
                    \stackrel{b_k}{\tto}\  x_k\ =:\ {\rm end}(\omega)\,,
 $$
with no descents, that is $b_1\leq b_2\leq\cdots\leq b_k$.
The resulting quasi-symmetric function $K_P$ was studied in~\cite{BS99b},
where (with some effort), it was shown to give a Hopf morphism from a
reduced incidence Hopf algebra to $\Q$.
We may likewise have edge-labelled graphs, and define descent Pieri operators
in that context.

To a subset $\{j_1<j_2<\cdots<j_k\}$ of $[n-1]$,
we associate the composition $(j_1,j_2-j_1,\ldots,n-j_k)$.
Given a saturated chain $\omega$ in $P$ with labels
$b_1,b_2,\ldots,b_n$, let $D(\omega)$ be the {\sl descent composition} of
$\omega$, that is the composition associated to the descent set
$\{i\mid b_i>b_{i+1}\}$ of $\omega$.
Then (Equation 4 of~\cite{BS99b}) we have
  \begin{equation}\label{eq:K-chains}
    K_{[x,y]}\ =\ \sum F_{D(\omega)}\,,
  \end{equation}
where the sum is over all saturated chains $\omega$ in the interval $[x,y]$,
and $F_\alpha$ is the complete (or fundamental) quasi-symmetric function. 

If we label a cover $\mu\lessdot\lambda$ in Young's lattice consistently by
either the column or content of the box in $\lambda/\mu$, then the descent
Pieri operator coincides with the Pieri operator of
Example~\ref{Ex:SkewSchur}.
\end{defn}

\begin{example}
$k$-Bruhat order and skew Schubert functions. \quad\rm
The Pieri-type formula for the classical flag manifold~\cite{LS_82,So96}
suggests a symmetric Pieri operator on a suborder of the Bruhat order on the
symmetric group, which encodes the structure of the cohomology of the flag
manifold as a module over the ring of symmetric polynomials.
Let ${\mathcal S}_n$ denote the symmetric group on $n$ elements and let
$\ell(w)$ be the length of a permutation $w$ in this Coxeter group.

We define the $k$-Bruhat order $<_k$ by its covers.
Given permutations $u,w\in{\mathcal S}_n$, we say that
$u\lessdot_k w$ if $\ell(u)+1=\ell(w)$ and
$u^{-1}w=(i,j)$, where $(i,j)$ is a reflection with $i\leq k<j$.
When $u\lessdot_k w$, we write $wu^{-1}=(a,b)$ with $a<b$ and label the
cover $u\lessdot_k w$ in the $k$-Bruhat order with the integer $b$.

The descent Pieri operators on this labelled poset are symmetric as
$\h_m$  models the action of the Schur polynomial $h_m(x_1,\ldots ,x_k)$ on
the basis of Schubert classes (indexed by ${\mathcal S}_n$) in the
cohomology of the flag manifold $SL(n,\mathbb{C})/B$.
We also have
 $$
  K_{[u,w]}\ =\ \sum_\lambda c^w_{u,(\lambda,k)}\, S_\lambda
 $$
where $c^w_{u,(\lambda,k)}$ is the coefficient of the Schubert polynomial
${\mathfrak S}_w$ in the product\break
${\mathfrak S}_u\cdot S_\lambda(x_1,\ldots,x_k)$.
This is the skew Schubert function $S_{wu^{-1}}$ of~\cite{BS_skew}.
Geometry shows these coefficients $c^w_{u,(\lambda,k)}$ are non-negative.
It is an important open problem to give a combinatorial or algebraic proof
of this fact.
\end{example}

\begin{example}\label{Ex:Stanley}
The weak order on ${\mathcal S}_n$ and Stanley symmetric functions.
\quad\rm
The {\sl weak order} on the symmetric group ${\mathcal S}_n$
is the labelled poset whose covers are $w\lessdot w(i,i+1)$,
with label $i$ if $\ell(w)+1=\ell(w(i,i+1))$.
In~\cite{BS99b}, it is shown that the descent Pieri operators on this
labelled poset are symmetric and $\F_{[u,w]}$ is the Stanley symmetric
function or stable Schubert polynomial $F_{wu^{-1}}$, introduced by Stanley
to study reduced decompositions of chains in the weak order on
${\mathcal S}_n$~\cite{Stanley84}.
\end{example}

\begin{example}
noncommutative Schur functions of Fomin and Greene.
\quad \rm
Fomin and Greene have a theory of {\sl combinatorial representations} of
certain noncommutative Schur functions~\cite{FG}.
These are a different noncommutative version of symmetric
functions than $NC$.
Using the Cauchy element in their algebra, they obtain
symmetric functions $F_{y/x}$ which include Schur functions, Stanley
symmetric functions, stable Grothendieck polynomials, and others.
A combinatorial representation gives rise to
an edge-labelled directed graph so that the
functions $F_{y/x}$ of Fomin-Greene are the functions $K_{[x,y]}$
coming from the descent Pieri operators on this structure.

Let $FG_n$ be the quotient of the free
associative algebra $\Z\Span{u_1,u_2,\ldots,u_n}$ by the two-sided ideal
generated by the following relations
  \begin{equation}\label{eq:FG-relations}
    \begin{array}{rclll}
      u_iu_ku_j&=&u_ku_iu_j,&i\leq j<k& |i-k|\geq 2\\
      u_ju_iu_k&=&u_ju_ku_i,&i<j\leq k& |i-k|\geq  2\\
      (u_i+u_{i+1})u_{i+1}u_i&=&u_{i+1}u_i(u_i+u_{i+1}).
    \end{array}
  \end{equation}

In $FG_n\times \mathbb{Z}[z_1,z_2,\ldots,z_m]$ define the noncommutative
Cauchy element to be
 $$
  \psi\ :=\ \prod_{i=1}^m\prod_{j=n}^1(1+z_iu_j).
 $$
Let $R$ be any set whose cardinality is at most countable, and let $\Z R$ be
the free abelian additive group with basis consisting of the elements of
$R$.
A representation  of $FG_n$ on $\Z R$ is {\sl combinatorial} if for
all $x\in R$, we have $x.u_i\in R\cup \{ 0\}$.
Given a combinatorial representation of $FG_n$ on $\Z R$ and $x,y\in R$, set
 $$
  F_{y/x}\ :=\ \Span{x.\psi,\ y}.
 $$

We define an edge-labelled directed multigraph $\mathfrak{R}$ with vertex set
$R$ for which  $F_{y/x}$ is the quasi-symmetric function coming from the
descent Pieri operator on that structure.
We construct   $\mathfrak{R}$   by drawing an edge
with label $-i$ from $x$ to $x.u_i$
 $$
  x  \stackrel{-i}{\ttto} x.u_i
 $$
if $x.u_i\neq 0$.
Considering the descent Pieri operators on $\mathfrak{R}$,
we have the following.

  \begin{thm}\label{thm:FG}
     For every $x,y\in R$, $F_{y/x}=\F _{[x,y]}(z_1,\ldots,z_m,0)$.
  \end{thm}

\begin{Remark}\rm
We identify the generators $z_i$ in $\Z[z_1,\ldots,z_m]$ with those in the
algebra $\mathcal{Q}sym$ generated by
the indeterminates $z_1,z_2,\ldots$.
\end{Remark}

Before we prove Theorem ~\ref{thm:FG}, we recall some results from \cite{FG}.
Define
 $$
  e_k({\bf u})\ :=\
   \sum_{i_1>i_2>\cdots>i_k}u_{i_1}u_{i_2}\cdots u_{i_k}
 $$
and for a partition $\lambda=(\lambda_1,\ldots,\lambda_m)$ set
$e_\lambda({\bf u}):=e_{\lambda_1}({\bf u})\cdots e_{\lambda_m}({\bf u})$.

\begin{prop}[Fomin-Greene]\label{prop:FG}
\mbox{ }
  \begin{enumerate}
   \item[(1)]
      For any positive integers $a,b$, we have
      $e_a({\bf u})e_b({\bf u})=e_b({\bf u})e_a({\bf u})$.
   \item[(2)]
      $\psi=\sum_\lambda m_\lambda(z)e_\lambda({\bf u})$.
  \end{enumerate}
Here, $m_\lambda(z):= m_\lambda(z_1,\ldots,z_m)$ is the monomial
symmetric polynomial.
\end{prop}

\noindent{\it Proof of Theorem~\ref{thm:FG}.}
Observe that for $x\in R$,
 $$
  x.e_k({\bf u})\ =\
  \sum_{\stackrel{\mbox{\scriptsize $x\stackrel{{-i_1}}{\ttto}
  \cdots\stackrel{{-i_k}}{\ttto}y$}}{i_1>\cdots>i_k}} y
   \ =\ x.\overline{h}_k.
 $$
{From} this, and Proposition ~\ref{prop:FG} it follows that the Pieri
operators are symmetric, that is $\h_a\h_b=\h_b\h_a$ for all $a,b\in \Z ^+$.
Hence, as in Example~\ref{ex:sym},
we have $x.S^\alpha=x.S^\beta$ whenever $\alpha$ and $\beta$ are two
compositions that
determine the same partition.

Then
  \begin{eqnarray*}
     F_{y/x}\ =\ \langle x.\psi ,y\rangle&=&
      \sum_\lambda m_\lambda(z)\langle x.e_\lambda(\textbf{u}),y\rangle\\
    &=&\sum_\lambda m_\lambda(z)\langle x.S^\lambda ,y\rangle \\
    &=&\langle x. \sum_\alpha M_\alpha(z )S^\alpha,y\rangle \\
    &=&\sum _r\langle x. \sum_{\alpha\vDash r} M_\alpha(z )S^\alpha
         ,y\rangle \\
    &=&\sum _r \F_{[x,y]^{(r)}}(z_1,\ldots ,z_m,0)
 \end{eqnarray*}
by Theorem ~\ref{GCauchy},
which by definition is equal to $ \F _{[x,y]}(z_1,\ldots,z_m,0)$.
\QED

\end{example}

%%%%%%%%%%%%%%%%%%%%%%%
\begin{example}
$P$-Partitions.  \quad\rm
 Let $P$ be a poset and consider any (vertex) labelling $\gamma\colon
P\rightarrow {\mathbb N}$   of $P$.
 A $(P,\gamma)$-{\sl partition} is an order preserving function
$f\colon P\rightarrow {\mathbb N}$
such that if $x<y$ and $\gamma(x)>\gamma(y)$, then $f(x)<f(y)$.
It is sufficient to check these conditions for covers $x\lessdot y$ in $P$.

Let ${\mathcal A}(P,\gamma)$ be the set of all
$(P,\gamma)$-partitions.
The {\sl weight enumerator}
$\Gamma(P,\gamma)$ of the labelled poset $(P,\gamma)$ is
 $$
   \Gamma(P,\gamma)\ :=\
    \sum_{f\in{\mathcal A}(P,\gamma)}\  \prod_{x\in P} z_{f(x)}\,.
 $$
This is obviously quasi-symmetric.

Properties of this weight enumerator are tied up with
Stanley's Fundamental
Theorem of $P$-partitions~\cite{Stanley72}.
Let ${\mathcal L}(P)$ be the set of all linear extensions of $P$.
A linear extension $w$ of $P$ lists
the elements of $P$ in order $w_1,w_2,\ldots,w_n$, with $w_i<w_j$ (in $P$)
implying $i<j$.
Here $n = |P|$.
Let $D(w,\gamma)$ be the descent composition of $n$ associated to the descent
set of the sequence of integers
$\gamma(w_1),\gamma(w_2),\ldots,\gamma(w_n)$.

For a linear ordering $w$ of $P$, the set ${\mathcal A}(w,\gamma)$ may be
identified with the set of all weakly increasing functions
$f\colon [n]\rightarrow{\mathbb N}$ where if $\gamma(w_i)>\gamma(w_{i+1})$
then $f(w_i)<f(w_{i+1})$.
Thus $\Gamma(w,\gamma)$ is Gessel's fundamental
quasi-symmetric function~\cite{Ges}
$F_{D(w,\gamma)}$.

The Fundamental Theorem of $P$-partitions notes that
 $$
  {\mathcal A}(P,\gamma)\ =\
  \coprod_{w\in{\mathcal L}(P)} {\mathcal A}(w,\gamma)\,.
 $$
This implies that
  \begin{equation}\label{eq:FundP}
     \Gamma(P,\gamma)\ =\ \sum_{w\in{\mathcal L}(P)} \Gamma(w,\gamma)\ =\
     \sum_{w\in{\mathcal L}(P)} F_{D(w,\gamma)}\,.
  \end{equation}
We  show that $\Gamma(P,\gamma)$
is given by descent Pieri operators on the (graded) poset ${\mathcal I}P$ of
lower order ideals of P with (edge) labelling induced from the vertex
labelling of $P$.
A subset $I\subset P$ is a {\sl lower order ideal} of $P$
if whenever $x\in I$ and $y<x$, then $y\in I$.
The set ${\mathcal I}P$ of lower order ideals of $P$ is ordered by inclusion.
We label a cover $I\subset\hspace{-7pt}\cdot\hspace{3pt} J$ in
${\mathcal I}P$ with $\gamma(x)$, where $x$ is the unique
element $x \in J\setminus I$.
Then ${\mathcal L}(P)$ is in bijection
with the maximal chains of ${\mathcal I}P$.
Using the descent Pieri operators for this structure,
Equation~\ref{eq:K-chains} shows that every maximal chain
of ${\mathcal I}P$ contributes the summand $F_{D(w,\gamma)}$
to $\F_{{\mathcal I}P}$ where $w$ is the linear extension
of that chain.
Thus
 $$
  \F_{{\mathcal I}P}\ =\ \Gamma(P,\gamma)\,.
 $$
The Hopf structure of $\H ({\mathcal I}P)$
was studied by Malvenuto in~\cite{Malv}.
\end{example}

\begin{example}
Quantum cohomology of Grassmannian, fusion coefficients, and the Hecke
algebra at a root of unity.
\quad\rm
Let $m,p$ be positive integers and let ${\mathcal C}_{p,m}$ be
the set of sequences $\alpha:0<\alpha_1<\cdots<\alpha_p$ which also satisfy
$\alpha_p-\alpha_1<m+p$.
We order this set of sequences by componentwise comparison to obtain a
ranked poset.
Given a cover $\alpha\lessdot\beta$, there is a unique index $i$ with
$\alpha_i+1=\beta_i$ and $\alpha_j=\beta_j$ for $i\neq j$.
We label such a cover with $\beta_i$.

The elements of the poset ${\mathcal C}_{m,p}$ may alternately be
described by pairs $(a,\lambda)$, where $a$ is a positive integer and
$\lambda$ is a partition with $\lambda_{p+1}=0$ and $\lambda_1\leq m$.
We obtain $(a,\lambda)$ from the sequence $\alpha$ by
  \begin{eqnarray*}
     \{\lambda_1+p,\ldots,\lambda_p+1\}&\equiv&
         \{\alpha_1,\ldots,\alpha_p\}\ \mbox{mod}\ (m+p), \\
      a\cdot(m+p)&=&\sum_{i=1}^p \alpha_i-\lambda_i-i\,.\rule{0pt}{15pt}
  \end{eqnarray*}
We may likewise pass from the indexing scheme $(a,\lambda)$ to sequences
$\alpha$, as this association is invertible (see~\cite{RRW96}).

For $x\in{\mathcal C}_{m,p}$ and $0<k\leq p$, consider the Pieri operator
  \begin{equation}\label{eq:qpf}
    x.\h_k\ :=\ \sum_{\omega} \mbox{end}(\omega),
  \end{equation}
where the sum is over all chains $\omega$ of length $k$ starting at $x$,
 $$
  \omega \ :\  x  \ \stackrel{b_1}{\tto}\
               x_1\ \stackrel{b_2}{\tto}\  \cdots
                    \stackrel{b_k}{\tto}\  x_k\ =:\ {\rm end}(\omega)\,,
 $$
with no descents, that is $b_1\leq b_2\leq\cdots\leq b_k$, and also
satisfying the restriction $b_k-b_1<m+p$.
Thus these operators $\h_k$ are not an instance of rank-selection or descent
Pieri operators as previously introduced.

These Pieri operators $\h_k$ are symmetric as they model the Pieri formula
in the quantum cohomology ring~\cite{Be97} of the Grassmannian
of $p$-planes in ${\mathbb C}^{m+p}$.
This commutative quantum cohomology ring has a basis $q^a\sigma_\lambda$ for
$(a,\lambda)\in{\mathcal C}_{m,p}$, and
 $$
  q^a \sigma_\lambda\cdot\sigma_k \ =\
   \sum_{(b,\mu)} q^b \sigma_\mu\,,
 $$
where the sum is over all indices $(b,\mu)$ appearing in the product
$(a,\lambda).\h_k$~(\ref{eq:qpf}), when it is written in terms of pairs.
Thus we have the following formula
 $$
   K_{[(b,\mu),(a,\lambda)]}\ =\
   \sum_\nu c^\lambda_{\mu,\nu}\, S_\nu\,,
 $$
where the sum is over all partitions $\nu$ of
$rk(a,\lambda)-rk(b,\mu)$ with $\nu_{p+1}=0$ and $m\geq\nu_1$.
Here, $c^\lambda_{\mu,\nu}$ is the quantum Littlewood-Richardson
coefficient~\cite{BCF}, the coefficient of $q^a\sigma_\lambda$ in the
product $q^b\sigma_\mu\cdot \sigma_\nu$.

These Pieri operators also model the fusion product in the Verlinde algebra
(see~\cite{BCF} for a discussion), and the Pieri formula in the
representation rings of Hecke algebras at roots of unity~\cite{GW90}.
Geometry and representation theory show that these coefficients
$c^\lambda_{\mu,\nu}$ are non-negative, but a combinatorial
proof of this fact is lacking.
\end{example}

%%%%%%%%%%%%%%%%%%%%%%%%%%%%%%%%%%%%%%
\section{Peak enumeration and Eulerian Pieri operators}

\begin{defn}\label{Ex:PPO}
Let $\omega$ be a labelled ordered chain, that is
 $$
  \omega \ :\  x_0  \ \stackrel{b_1}{\tto}\
               x_1\ \stackrel{b_2}{\tto}\  \cdots
                    \stackrel{b_k}{\tto}\  x_k\,.
 $$
We say that $\omega$ has a {\sl peak} at $i$ if
$b_{i-1}\leq b_{i}>b_{i+1}$.
Let $\Lambda(\omega)$ be the {\sl peak composition} of $\omega$, that is the
composition
of $k$ associated to the peak set $\{i | b_{i-1}\leq b_{i}>b_{i+1} \}$ of
$\omega$. Let $P$ be a labelled poset.
To enumerate chains in intervals $[x,y]$ of $P$ according to
their peaks, we use the peak enumerator
 $$
  \widetilde{K}_{[x,y]}\ :=\ \sum_\omega \theta_{\Lambda(\omega)}\,,
 $$
where the sum is over all saturated chains $\omega$ in the interval $[x,y]$.
We show this peak enumerator is the quasi-symmetric function
$K_{\delta[x,y]}$ associated to the descent Pieri operators on an enriched
structure $\delta P$ defined on the labelled poset $P$.

Given a labelled poset $P$, where (for simplicity)  we assume that the labels
$b_i$ are positive integers, we define $\delta P$, the {\sl doubling} of $P$,
to be the labelled directed graph with vertex set $P$, where
every edge $x\stackrel{b}{\tto}y$ of $P$ is doubled, but with one label the
negative of the original label, that is
 $$
  \begin{picture}(66,24)(0,3)
   \put(0,13){$x$}
   \put(55,13){$y\ .$}
   \put(18,16){$\ttto$}
   \put(18, 9){$\ttto$}
   \put(29,21){$\scriptstyle b$}
   \put(22,4){$\scriptstyle -b$}
  \end{picture}
 $$
Such a  poset whose Hasse diagram has
 multiple edges is called a {\sl r\'eseau}. The r\'eseau $\delta P$ is the
{\sl doubled r\'eseau}
of $P$.
To define descent Pieri operators on the r\'eseau $\delta P$, we say that
there is a descent at $i$ if consecutive
labels $b_i,b_{i+1}$ satisfy either $b_i>b_{i+1}$ or else
$b_i=b_{i+1}<0$. We then adjust the definitions of descent set and descent
composition
accordingly.

The following Theorem is a generalisation 
of~\cite[Theorem 3.6]{Stembridge_enriched}, as will becom apparent from
Example 7.5. 

  \begin{thm}\label{thm:peak-descent}
   Let $P$ be any labelled poset and $\delta P$ its doubled r\'eseau.
   Then the modified descent Pieri operators on $\delta P$ are Eulerian, and we
   have
   $$
    K_{\delta[x,y]}
    \ =\ \widetilde{K}_{[x,y]}\ =\ \sum c_{x,\alpha}^y\, \theta_\alpha\,.
   $$
   where the sum is only over peak compositions $\alpha$.
  \end{thm}

These combinatorial invariants
$c_{x,\alpha}^y$ of $\delta P$ enumerate the chains of $P$
whose peak sets have composition $\alpha$.

Before we prove Theorem~\ref{thm:peak-descent}, we make some definitions and
prove two auxiliary lemmas.
For a composition $\alpha$ of $n$, let $\alpha^+$ be the composition of $n+1$
obtained from $\alpha$ by increasing its last component by 1, and
$\alpha\cdot 1$ be the composition of $n+1$ obtained by appending a
component of size 1 to $\alpha$.
Define linear maps $\psi,\varphi\colon\Q_n\to \Q_{n+1}$ by
  \begin{eqnarray*}
    \psi(M_\beta)&:=& M_{\beta^+} + 2  M_{\beta\cdot 1}\,,\\
    \varphi(M_\beta)&:=&
            \delta_{1,\beta_l}M_{\beta^+}+2 M_{\beta\cdot 1}\,,
  \end{eqnarray*}
where $\beta_l$ is the last component of $\beta$ and $\delta_{1,\beta_l}$
is the
Kronecker delta function.
Using the relation $F_\beta=\sum_{\alpha\preccurlyeq\beta}M_\alpha$ between
the two bases of $\Q$, we see that
 $$
  \psi(F_\beta)\ =\ F_{\beta^+} + F_{\beta\cdot 1}\,.
 $$

  \begin{lem}\label{lem:thetamap}
    $\psi(\theta_\alpha)=\theta_{\alpha^+}$ \  and \
    $\varphi(\theta_\alpha)=\theta_{\alpha\cdot 1}$.
  \end{lem}

\proof
The function $\theta_{\alpha\cdot 1}$ is the sum of terms
$2^{k(\beta)}M_\beta$ for each $\beta$ satisfying
$\beta^\ast\preccurlyeq\alpha\cdot 1$.
Suppose $\beta^*\preccurlyeq \alpha\cdot 1$.
If $\beta=\gamma\cdot 1$, then $\beta^*=\gamma^*\cdot 1$ and we have
$\gamma^*\preccurlyeq \alpha$.
Conversely, if $\gamma^*\preccurlyeq \alpha$, then $\beta:=\gamma\cdot 1$
satisfies $\beta^*\preccurlyeq \alpha\cdot 1$.
Thus every summand $2^{k(\gamma)}M_\gamma$ of $\theta_\alpha$ contributes a
summand $2\cdot 2^{k(\gamma)}M_{\gamma\cdot 1}$ to $\theta_{\alpha\cdot 1}$.

The other summands $\beta$
have  $\beta=\gamma\cdot\beta_l$ with $\beta_l>1$.
Then $\beta^*=\gamma^*\cdot 1\cdot(\beta_l-1)$.
If $\beta^*\preccurlyeq \alpha\cdot 1$, then we must have $\beta_l=2$, so that
$\beta=\gamma\cdot 2$.
Then  $\beta^*=\gamma^*\cdot 1\cdot 1\preccurlyeq\alpha\cdot 1$, which
implies that $(\gamma\cdot 1)^*\preccurlyeq \alpha$.
Conversely, if $(\gamma\cdot 1)^*\preccurlyeq \alpha$, then
$(\gamma\cdot 2)^*\preccurlyeq \alpha\cdot 1$.
Thus every summand $2^{k(\gamma\cdot 1)}M_{\gamma\cdot 1}$ of
$\theta_\alpha$ contributes a summand $2^{k(\gamma\cdot 2)}M_{\gamma\cdot 2}$
to $\theta_{\alpha\cdot 1}$.

This shows that $\theta_{\alpha\cdot 1}=\varphi(\theta_\alpha)$.
The arguments for $\theta_{\alpha^+}$ are similar, but simpler.
\QED

The key lemma relating the peak enumerator on $P$ and the modified descent
Pieri operators on the r\'eseau $\delta P$ concerns the case when $P$ is a
chain.

  \begin{lem}\label{lem:peak_chains}
   Suppose $\omega$ is a chain.
   Then $K_{\delta\omega}=\theta_{\Lambda(\omega)}$.
  \end{lem}

\proof
We prove this by induction on the length of the chain $\omega$.
The initial cases are easy calculations.
Let $b_1,\ldots,b_k$ be the word of $\omega$, and set $u$ to be the
truncation of $\omega$ at the penultimate cover, so that $b_1,\ldots,b_{k-1}$
is the word of $u$.

Consider first the case where $b_{k-1}\leq b_k$.
Then every chain $\gamma$ in $\delta u$ gives two chains $\gamma.b_k$ and
$\gamma.\overline{b_k}$ in $\delta \omega$.
Since $D(\gamma.b_k)=D(\gamma)^+$ and
$D(\gamma.\overline{b}_k)=D(\gamma)\cdot 1$, we see that
$K_{\delta\omega}=\psi(K_{\delta u})$.
Similarly, if $b_{k-2}>b_{k-1}>b_k$, then considering the last
three labels of a chain in $\delta \omega$ show that
$K_{\delta\omega}=\psi(K_{\delta u})$.
In both cases, $\Lambda(\omega)=\Lambda(u)^+$ (as the peak sets are the
same),  and the lemma follows by Lemma~\ref{lem:thetamap}.

Now suppose $b_{k-2}\leq b_{k-1}>b_k$.
Let $v$ be the truncation of $\omega$ at the $(k-2)$th position.
Let $\gamma$ be a chain of $\delta v$ with descent composition $\alpha$.
Then $\gamma$ has 4 extensions to chains in $\delta\omega$, and 2 have
descent composition $\alpha^+\cdot 1$ and 2 have descent composition
$\alpha\cdot 2$.
Thus if we define
$\phi\colon F_\alpha\mapsto 2F_{\alpha^+\cdot 1}+2F_{\alpha\cdot 2}$, then
$\phi(K_{\delta v})=K_{\delta\omega}$.
A straightforward calculation shows
$\phi(M_\beta)=2M_{\beta^+\cdot 1}+2 M_{\beta\cdot 2}
               +4M_{\beta\cdot 1\cdot 1}$,
which is $\varphi(\psi(M_\beta))$.
Thus
$K_{\delta\omega}=\varphi(\psi(K_{\delta v}))=\varphi(K_{\delta u})$.
Since $\omega$ has a peak at $n-1$, we have
$\Lambda(\omega)=\Lambda(u)\cdot 1$, and so this case follows by
Lemma~\ref{lem:thetamap}.
\QED

\noindent{\it Proof of Theorem~\ref{thm:peak-descent}.}
Given an interval $[x,y]$ in a poset or r\'eseau, let
$\mbox{\it ch}[x,y]$ be the set of saturated chains in $[x,y]$.
Let $K$ be the quasi-symmetric function given by the descent Pieri operators
on the r\'eseau $\delta P$.
Let $x\leq y$ in $P$.
Given a chain $\omega\in\mbox{\it ch}\,\delta[x,y]$, we obtain a chain
$|\omega|\in\mbox{\it ch}[x,y]$ by replacing each cover in $\delta[x,y]$
with a negative integer label by the corresponding cover in $[x,y]$ whose
label is positive.
Then, by Equation~\ref{eq:K-chains}, we have
  \begin{eqnarray*}
  \F_{\delta[x,y]}
    &=& \sum_{\omega\in\mbox{\scriptsize\it ch}\,\delta[x,y]}F_{D(\omega)}\\
    &=& \sum_{\beta\in\mbox{\scriptsize\it ch}[x,y]}
          \  \sum_{\omega : |\omega|=\beta} F_{D(\omega)}\\
    &=& \sum_{\beta\in\mbox{\scriptsize\it ch}[x,y]} \F_{\delta \beta}
     \quad =\quad
     \sum_{\beta\in\mbox{\scriptsize\it ch}[x,y]} \theta_{\Lambda(\beta)}
     \quad =\quad  \widetilde{K}_{[x,y]}\,.\qquad\QED
  \end{eqnarray*}
\end{defn}

\begin{example}
Enriched $P$-partitions. \quad\rm
Stembridge enriches the theory of $P$-parti\-tions~\cite{Stembridge_enriched}
giving a new class of quasi-symmetric generating functions.
Let $(P,\gamma)$ be a labelled poset and let
$\PP=\{\overline{1},1,\overline{2},2,\overline{3},3,\ldots\}$ be two copies
of the positive integers ordered as follows:
$\overline{1}<1<\overline{2}<2<\overline{3}<3<\cdots$.
An {\sl enriched $(P,\gamma)$-partition}
is an order-preserving map $f\colon P\to\PP$ such that for $x<y$ in $P$ and
$k\in\Z^+$
  \begin{itemize}
   \item[-] if $f(x)=f(y)=\overline{k}$, then $\gamma(x)<\gamma(y)$,
   \item[-] if $f(x)=f(y)={k}$, then $\gamma(x)>\gamma(y)$.
  \end{itemize}

\noindent
Let ${\mathcal E}(P,\gamma)$ be the set of all enriched
$(P,\gamma)$-partitions and  define the weight enumerator
 $$
  \Delta(P,\gamma)\ =\
   \sum_{f\in{\mathcal E}(P,\gamma)}\  \prod_{x\in P}  z_{f(x)}
 $$
where $z_{\overline k}=z_k$ for all positive integers $k$.
The analogue of Equation~\ref{eq:FundP} for enriched $P$-partitions
 is
 $$
   \Delta(P,\gamma)\ =\
   \sum_{w\in{\mathcal L}(P)} \Delta(w,\gamma).
 $$
We thus need to characterise the quasi-symmetric function corresponding to
a linear extension $(w,\gamma)$ of $(P,\gamma)$.
A {\it peak} of the linear extension $(w,\gamma)$ is an index
$i$ with $1 < i < |w|$ where
$\gamma(w_{i-1})<\gamma(w_i)>\gamma(w_{i+1})$.
Stembridge shows that
$$
  \Delta(w,\gamma)\ =\ \theta_{\Lambda(w,\gamma)}
$$
where
$\Lambda(w,\gamma)$ is the peak composition associated to the
peak set of the linear extension $(w,\gamma)$.
We can then generalise the construction we have for $P$-partitions.
This time we  proceed as in Definition~\ref{Ex:PPO}
and
consider the descent Pieri operators on the doubled r\'eseau
$\delta{\mathcal I}P$.
By Lemma~\ref{lem:peak_chains}, every maximal chain of
${\mathcal I}P$ contributes exactly
$\Delta(w,\gamma)$ to  $\F_{\delta{\mathcal I}P}$,
where $w$ is the linear extension of
${\mathcal L}(P)$ corresponding to that chain.
This shows that
 $$
  \F_{\delta{\mathcal I}P}\ =\ \Delta(P,\gamma).
 $$
\end{example}

\begin{example}\label{Ex:iso}
Isotropic Pieri formula.  \quad\rm
The Pieri-type formulas for the flag manifolds $SO(2n+1,{\mathbb C})$
and $Sp(2n,{\mathbb C})$ of~\cite{BS_lag-pieri} each give
symmetric Eulerian Pieri operators.
These are defined on enrichments of the same subposet of the Bruhat
order on the group ${\mathcal B}_n$ of signed permutations.
For an integer $i$, let $\overline{\imath}$ denote $-i$.

We regard ${\mathcal B}_n$ as a subgroup of the group of permutations
on $\{\overline{n},\ldots,\overline{2},\overline{1},1,\ldots,n\}$.
Let $\ell$ be the length function on the Coxeter group ${\mathcal B}_n$.
The {\sl $0$-Bruhat order} $<_0$ on ${\mathcal B}_n$ is the labelled poset
${\mathcal B}^0_n$
with covers $u\lessdot_0w$ if
$\ell(u)+1=\ell(w)$ and $u^{-1}w$ is a reflection with either the form
$(\overline{\imath},i)$ or the form
$(\overline{\imath},j)(\overline{\jmath},i)$ for some $0<i,j$.
When $u\lessdot_0w$, either $wu^{-1}=(\overline{\beta},\beta)$ for some
$0<\beta$ or else
$wu^{-1}=(\overline{\beta},\overline{\alpha})(\alpha,\beta)$ for some
$0<\alpha<\beta\leq n$.
We label a such a cover with the (positive) integer $\beta$.

Consider the Eulerian descent Pieri operators on the doubled
r\'eseau $\delta{\mathcal B}_n^0$.
These operators are symmetric, as $\frac{1}{2}\h_k$ models the action
of the Schur P-polynomial $p_k$ on the basis of
Schubert classes (indexed by ${\mathcal B}_n$) in the cohomology of the flag
manifold  $SO(2n+1,{\mathbb C})/B$~\cite{BS_lag-pieri}.
(This is because there are twice as many increasing chains in a doubled
interval $\delta[x,y]$ as peakless chains in the interval $[x,y]$, and the
coefficient of $y$ in $x.p_k$ is this number of peakless chains.)

We modify this descent action of $NC$ on $\delta{\mathcal B}_n^0$
by identifying $h_k$ with $\frac{1}{2}\h_k$, which is still integral.
These new Pieri operators are symmetric, as they model the action
of $p_k$, {\it and} they are Eulerian, as $2h_k$ satisfies the Euler
relations~\ref{eq:Euler}.
In exact analogy to how the Skew Schubert functions are shown
in~\cite{BS_skew} to be the generating functions for the coefficients
$c^w_{u,(\lambda,k)}$, given by descent Pieri operators, we have
the following formula
 $$
  K_{[u,w]}\ =\ \sum_\lambda b^w_{u,\lambda}\, Q_\lambda\,,
 $$
where the sum is over all strict partitions $\lambda$ of $\ell(w)-\ell(u)$.
Here $b^w_{u,\lambda}$ is the coefficient of the Schubert class
${\mathfrak B}_w$ in the product ${\mathfrak B}_u\cdot P_\lambda$,
and $P_\lambda,Q_\lambda$ are Schur P- and Q-polynomials, which form dual
bases for the self dual symmetric Hopf algebra $\Pi_{\mathbb Q}\cap \Lambda$.
The polynomials $Q_\lambda$ appear as
$\sum_\lambda P_\lambda\otimes Q_\lambda$ is the Cauchy element of
$\Pi_{\mathbb Q}\cap\Lambda$.

For the symplectic flag manifold, we modify the r\'eseau
$\delta{\mathcal B}_n^0$ by erasing the negative edge in a cover
%%% A picture in the text
 %% This adds enough space above and below the line
   \rule[-10pt]{0mm}{8pt}
   $\begin{picture}(47,14)(3,0)
     \put(2,0){$u$}
     \put(40,0){$w$}
     \put(15,4){\vector(1,0){20}}
     \put(15,1){\vector(1,0){20}}
     \put(22,8){$\scriptstyle \beta$}
     \put(15,-7){$\scriptstyle -\beta$}
    \end{picture}
   $
when $wu^{-1}=(\overline{\beta},\beta)$.
Write ${\mathcal L}{\mathcal B}^0_n$ for the resulting r\'eseau.
It is a slight modification of the 0-Bruhat r\'eseau of~\cite{BS_lag-pieri},
and may be used in its place for the combinatorics therin.
Let $\{\h_k\}$ be the descent Pieri operators on
${\mathcal L}{\mathcal B}^0_n$.
This family of Pieri operators is symmetric and Eulerian, as $\h_k$
models the action of the Schur Q-polynomial
$q_k$ on the Schubert basis of the cohomology of the flag manifold
$Sp(2n,{\mathbb C})$, and the Schur Q-polynomials $q_k$ satisfy the Euler
relations.
As before, we have the following formula
 $$
  K_{[u,w]}\ =\ \sum_\lambda c^w_{u,\lambda}\, P_\lambda\,,
 $$
where the sum is over all strict partitions $\lambda$ of $\ell(w)-\ell(u)$.
Here $c^w_{u,\lambda}$ is the coefficient of the Schubert class
${\mathfrak C}_w$ in the product ${\mathfrak C}_u\cdot Q_\lambda$.

Since every chain in an interval of ${\mathcal B}^0_n$ has the same number
of covers of the form $(\overline{\beta},\beta)$---these count the number
$s(wu^{-1})$ of sign changes between $u$ and $w$, we have
 $$
   \F_{\delta[u,w]}\ =\ 2^{-s(wu^{-1})} \F_{{\mathcal L}[u,w]}\,.
 $$
Geometry shows these coefficients $b^w_{u,\lambda}$ and $c^w_{u,\lambda}$
are non-negative.
It is an important open problem to give a combinatorial or algebraic proof
of this fact.
\end{example}

\begin{example}
Stanley symmetric functions of types $B$, $C$, and $D$.  \quad\rm
In \cite{BilHaim}, Billey and Haiman describe the Stanley symmetric
functions of types $B$ and $D$ in terms of peaks of reduced words of elements
in the corresponding Coxeter groups.

For ${\mathcal B}_n$, the simple transpositions are
$s_0,s_1,\ldots,s_{n-1}$, where $s_0=(\overline{1},1)$ and if $i>0$, then
$s_i=(\overline{i+1},\overline{\imath})(i,i+1)$.
The weak order on ${\mathcal B}_n$ is the labelled poset whose covers are
$w\lessdot ws_i$ with label $i+1$ if $\ell(w)+1=\ell(w s_i)$.
A reduced word ${\bf a}$ for $w$ is sequence of labels of a chain in
${\mathcal B}_n$ from the identity $e$ to $w$.
Billey and Haiman define the Stanley symmetric function of
type $B$ to be
 $$
  F^B_w\ :=\ \sum_{{\bf a}\in R(w)} \theta_{\Lambda({\bf a})}\,,
 $$
where $R(w)$ is the set of reduced words for $w$ and $\Lambda({\bf a})$
is the peak composition of the reduced word ${\bf a}$.
By Theorem~\ref{thm:peak-descent}, $F^B_w$ is the function $K_{\delta[e,w]}$
obtained from the Eulerian descent operators on the doubled r\'eseau
$\delta{\mathcal B}_n$.
Billey and Haiman establish the formula
 $$
  F^B_w\ =\ \sum_{\lambda} f^w_\lambda\, Q_\lambda\,,
 $$
where the sum is over all strict partitions $\lambda$ of $\ell(w)$,
and $f^w_\lambda$ counts the reduced words that satisfy a condition
imposed by the partition $\lambda$ (coming from the shifted Edelmann-Greene
correspondence~\cite{Ha92}).
Thus the Eulerian descent Pieri operators on $\delta{\mathcal B}_n$ are also
symmetric.

While Billey and Haiman do not define Stanley functions of type
$C$, one reasonably sets $F^C_w:=2^{-s(w)}F^B_w$, where $s(w)$ is the number
of sign changes in the permutation $w$.
This is just the number of $s_0$'s appearing in any reduced word of $w$.
Let the r\'eseau ${\mathcal L}{\mathcal B}_n$ be the
modification of the doubled r\'eseau $\delta{\mathcal B}_n$ where we erase
the edge with negative label $\overline{1}$ for covers $w\lessdot ws_0$.
Then every chain in an interval ${\mathcal L}[e,w]$ of
${\mathcal L}{\mathcal B}_n$ gives rise to $2^{s(w)}$ chains in
$\delta[x,y]$, each with the same descents as the original chain.
Then Equation~\ref{eq:K-chains} and Theorem~\ref{thm:peak-descent}
show that
 $$
  \F_{{\mathcal L}[x,y]}\ =\ 2^{-s(w)}\F_{\delta[x,y]}\ =\
   2^{-s(w)}F^B_w\ =\ F^C_w\,.
 $$
This shows these descent Pieri operators are Eulerian and symmetric.\bigskip

The Coxeter group ${\mathcal D}_n$ has simple reflections
$s_1,s_{\hat{1}},s_2,\ldots,s_{n-1}$.
The weak order on ${\mathcal D}_n$ is the labelled poset with cover
$w\lessdot ws_i$ labelled by $i$ if $\ell(w)+1=\ell(ws_i)$.
Here, we set $\hat{1}<1$.
A reduced word ${\bf a}$ for $w$ as before is a chain in ${\mathcal D}_n$
from $e$ to $w$.
Since $s_1$ and $s_{\hat{1}}$ commute, there are no occurences of $1\hat{1}1$
or $\hat{1}1\hat{1}$ in a reduced word, so changing all occurrences of
$\hat{1}$ to $1$ does not change the peaks in a reduced word, and the type $D$
Stanley symmetric functions of Billey and Haiman satisfy
$$
  F_w^D\ =\
   \sum_{{\bf a}\in R(w)} 2^{-o({\bf a})}\theta_{\Lambda({\bf a})}\,,
$$
where $o({\bf a})$ counts the number of occurrences of $1$ and $\hat{1}$
in the reduced word ${\bf a}$.
Let $\delta{\mathcal D}_n$ and ${\mathcal L}{\mathcal D}_n$ be the doubled
r\'eseau and its modification, erasing all edges with (negative) labels
$-1$ and $-\hat{1}$.

\begin{thm}
$F^D_w=\F_{{\mathcal L}[e,w]}$.
\end{thm}

\proof
By  Theorem~\ref{thm:peak-descent}, we have
$$
  F^D_w\ =\ \sum_{{\bf a}\in R(w)} 2^{-o({\bf a})}\F_{\delta{\bf a}}\,.
$$
The theorem follows from Equation~\ref{eq:K-chains} and
the following 1 to $2^{o({\bf a})}$ map from chains in
${\mathcal L}{\bf a}$ to chains in $\delta{\bf a}$,
which preserves descents.
When there are no subwords $1\hat{1}$ or $\hat{1}1$ in a chain in
${\mathcal L}{\bf a}$, simply
make all possible substitutions of negative and positive labels for
each occurrence of $1$ and $\hat{1}$.
If however, there is a subword $1\hat{1}$, then there is another chain
differing from the first only in that subword (having $\hat{1}1$ instead),
and the map uses the substitutions in both chains
  \begin{eqnarray*}
    \hat{1}1&\longmapsto&
       \hat{1}1, \ \overline{\hat{1}}1,\ \overline{1}\hat{1},\
       \overline{1}\overline{\hat{1}}\\
      1\hat{1}&\longmapsto&
       1\hat{1}, \ 1\overline{\hat{1}},\ \hat{1}\overline{1},\
       \overline{\hat{1}}\overline{1}\ .\qquad\QED
  \end{eqnarray*}
Lastly, we remark that these descent operators on
${\mathcal L}{\mathcal D}_n$ are Eulerian and symmetric, as Billey and
Haiman give a formula
$$
  F^D_w\ =\ \sum_{\lambda} e^\lambda_w\,Q_\lambda\,,
$$
where $e^\lambda_w$ is a rational number that counts certain weighted
reduced words.
\end{example}

%%%%%%%%%%%%%%%%%%%%%%%%%%%%%
\appendix
\section{Hopf algebras}

 A Hopf algebra is an algebra whose linear dual is also an algebra, with
some compatibility conditions.
They are important in representation theory (and in this paper)
because they act on tensor products of their representations.
The usefulness of Hopf algebras in combinatorics is apparent from the
ubiquity of their applications.
In this section, we summarize the basic notions of Hopf algebras.

A $\Z$-module $\H$ is a {\sl coalgebra} if there are two maps
$\Delta:\H\rightarrow\H\otimes\H$ (coproduct) and $\epsilon:\H\rightarrow
\Z$ (counit or
augmentation) such that the following diagrams commute
 $$
  \begin{picture}(140,58)
   \put(13,45){$\H$}	\put(97,45){$\H\otimes \H$}
   \put(0,0){$\H\otimes \H$}	\put(85, 0){$\H\otimes \H\otimes \H$}
   \put(53,51){\scriptsize $\Delta$}	\put( 8,25){\scriptsize $\Delta$}
   \put(47.5, 6){\scriptsize $\Delta\otimes 1$}
   \put(120,25){\scriptsize$1\otimes\Delta$}
   \put(30,48){\vector(1,0){62}} \put(40, 3){\vector(1,0){40}}
   \put(18,40){\vector(0,-1){27}}
   \put(116,40){\vector(0,-1){27}}
  \end{picture}
 \qquad\qquad
  \begin{picture}(107,58)
   \put(18,45){$\H$}	\put(72,45){$\H\otimes\H$}
   \put(6,0){$\H\otimes \H$}	\put(85,0){$\H$}
   \put(43,51){\scriptsize $\epsilon\otimes 1$}
   \put( 0,25){\scriptsize $1\otimes\epsilon$}
   \put(60, 6){\scriptsize$\Delta$}	\put(94,25){\scriptsize$\Delta$}
   \put(70,48){\vector(-1,0){40}} \put(77, 3){\vector(-1,0){32}}
   \put(23,13){\vector(0,1){27}}
   \put(90,13){\vector(0,1){27}} \put(80,10){\vector(-3,2){50}}
   \put(51,31){\scriptsize $1$}
  \end{picture}\qquad,
 $$
where $1$ is the identity map on $\H$.

\begin{rem}\rm The first of these diagrams is the coassociativity property,
which is the statement
that the dual of $\Delta$ defines an associative product on the linear dual
of $\H$, and the
second asserts this linear dual has a unit, induced by the dual of $\epsilon$.
\end{rem}

 If $\H$ is also an algebra, then it is a {\sl bialgebra} if
$\Delta,\epsilon$ are algebra morphisms.
While some authors call this structure a Hopf algebra, we define a
{\sl Hopf algebra} to be a bialgebra with a map
$s:\H\rightarrow \H$ (coinverse or antipode) such that the following diagram
commutes.
$$
  \begin{picture}(275,81)(-10,0)
  \put(-10,32){$\H$} \put(120,32){$\Z$}	\put(252,32){$\H$}
  \put( 50,67){$\H\otimes \H$}
  \put( 50, 0){$\H\otimes \H$} \put(165,67){$\H\otimes \H$}
  \put(165, 0){$\H\otimes \H$}
  \put(5,37){\vector(1,0){105}}
  \put(140,37){\vector(1,0){105}}
  \put(90,5){\vector(1,0){70}} \put(90,72){\vector(1,0){70}}
  \put(5,32){\vector(2,-1){40}}
  \put(205,62){\vector(2,-1){40}} \put(5,42){\vector(2, 1){40}}
  \put(205,12){\vector(2, 1){40}}
  \put(18,14){\scriptsize$\Delta$}\put(18,54){\scriptsize$\Delta$}
  \put(233,14){\scriptsize$\mu$}
  \put(233,54){\scriptsize$\mu$}
  \put(115,8){\scriptsize$1\otimes s$}
  \put(115,75){\scriptsize$s\otimes 1$}
  \put(65,40){\scriptsize$\epsilon$}
  \put(180,40){\scriptsize$u$}
  \end{picture}\qquad.
$$
Here
$\mu:\H\otimes\H\rightarrow \H$ is the
map induced by the
multiplication of $\H$ and $u:\Z\rightarrow \H$ is the map induced
by mapping $1$ to the unit of $\H$.
 The above diagram implies that $s$ is an algebra antihomomorphism, i.e.
$s(hh')=s(h')s(h)$ for all $h,h'\in \H$.

The existence of an antipode $s$ may seem to be a strong restriction on a
bialgebra, however as we will see,  it is no restriction for graded
bialgebras.
A {\sl graded bialgebra} is a  graded algebra $\H=\bigotimes_n\H_n$ where
$\Delta$ is graded and $\H_0=\Z$.
Given $x\in \H_n$, the $n$th graded component, we have
 $$
  \Delta(x)\ =\ x\otimes 1\ +\ \sum_{i=1}^n y_i\otimes z_{n-i},
 $$
where $y_i$ and $z_i$ have degree $i$.
The first term is always present due to the counit diagram. With this in
mind, Ehrenborg proved the following.

\begin{prop}[Lemma 2.1~\cite{Ehrenborg}]\label{GradedHopf}
 Given a graded bialgebra $\H$ there is a unique Hopf algebra with antipode
 $s$ defined recursively by $s(1)=1$, and for $x\in\H_n$, $n\geq 1$,
  $$
    s(x)\ =\ -\sum_{i=1}^n s(y_i)\cdot z_i.
  $$
\end{prop}

 Lastly, we remark on the useful Sweedler notation, which is an elegant
solution to the following quandary.
Given $h\in\H$, how do you efficiently represent $\Delta h$ as an element
of $\H\otimes\H$?
Carefully indexing this element would confuse even the writer.
Sweedler notation sidesteps this by omitting the indices of summation
entirely,
 $$
   \Delta h\ =\ \sum h_1\otimes h_2.
 $$
It is this notation that is normally used when dealing with Hopf algebras.

\end{document}